\documentclass{statsoc}

\usepackage{amsfonts}
\usepackage{mathrsfs}
\usepackage{graphicx}
\usepackage{amsmath}
\usepackage{verbatim}
\usepackage{natbib}
\usepackage{longtable}

\newtheorem{theorem}{Theorem}[section]
\newtheorem{lemma}[theorem]{Lemma}
\newtheorem{definition}[theorem]{Definition}

\newtheorem{example}[theorem]{Example}
\newtheorem{assumption}{Assumption}

\newtheorem{remark}[theorem]{Remark}

\title[Statistical analysis of the DH problem]{Statistical analysis of the
Diffie-Hellman key exchange protocol in a finite group\thanks{The
authors wish to thank Dr. Marco Lenci who suggested us the use of
the entropy function as a quantifier for random information.}}

\author[I. Florescu, A. Myasnikov and A. Mahalanobis]{Ionu\c{t} Florescu}
\address{Dept. of Mathematical Sciences, Stevens Institute of
Technology, Hoboken NJ 07030, USA} \email{ifloresc@stevens.edu}

\author[I. Florescu, A. Myasnikov and A. Mahalanobis]{Alexey Myasnikov}
\address{Dept. of Mathematical Sciences, Stevens Institute of
Technology, Hoboken NJ 07030, USA} \email{amyasnik@stevens.edu}

\author[I. Florescu, A. Myasnikov and A. Mahalanobis]{Ayan Mahalanobis}
\address{Dept. of Mathematical Sciences, Stevens Institute of
Technology, Hoboken NJ 07030, USA} \email{amahalan@stevens.edu}

\begin{document}

\input amssym.def \relax

\begin{abstract}
This paper presents a novel methodology to test the security of the
Diffie-Hellman public key exchange protocol. The security of many
cryptographic schemes rely on the hardness of this problem. We are
presenting a purely statistical test to compare this problem in
different groups. We are using groups included in $\mathbb Z_p$ with
$p$ prime as a major example, however the methods presented are not
restricted to these groups. The presentation of the results is
primarily intended to introduce novel applications of statistical
methodologies to the area of mathematical cryptography. As such we
will emphasize the cryptographical aspects of the work more than the
statistical notions.

\end{abstract}

\section{Introduction.}

Informally, through a key exchange protocol, two parties $\mathcal
A$ and $\mathcal B$ agree on a common key $K_{\mathcal A, \mathcal
B}$ pooled from a set $\mathcal S$ while communicating over an
insecure channel. Once the key is established, any further
information shared between the parties is encoded, transmitted and
decoded using the key $K_{\mathcal A, \mathcal B}$. The protocol is
secure if any third party $\mathcal C$ with access to the initial
communication between $\mathcal A$ and $\mathcal B$ cannot tell
apart $K_{\mathcal A, \mathcal B}$ from any other value in the set
$\mathcal S$. This guarantees that it is computationally unfeasible
for an outside adversary to gain ``any'' partial information on
$K_{\mathcal A, \mathcal B}$.

The Diffie-Hellman key exchange protocol \cite{DH} is a primary
example of a public key exchange protocol. In its most basic form,
the protocol chooses a finite cyclic group $(G,\cdot)$ of order $N$,
with generator $g$, where $\cdot$ denotes the group operation. In
what follows we chose the multiplicative operation to denote the
operation in the group, and thus the group $G$ is generated by the
powers of $g$ (i.e., $G=\{g^0, g^1,\dots,g^{N-1}\}$), symbolically
$G=<g>$. Note that $G$, $g$ and $N$ are public information.

The participants in the information transfer $\mathcal A$ and
$\mathcal B$ each randomly chooses an integer $a \in [1,N]$ and $b
\in [1,N]$ independently. Then $\mathcal A$ computes $g^a$,
$\mathcal B$ computes $g^b$ and exchange these elements of $G$ over
an insecure channel. Since each of $\mathcal{A}$ and $\mathcal{B}$
knows their respective $a$ and $b$ they can compute $g^{ab}$, which
or a publicly known derivation $K_{\mathcal{A,B}}$ of  that becomes
the public key.

Any method of converting $g^{ab}$ to $K_{\mathcal{A,B}}$ is publicly
known, and the security of the key $K_{\mathcal{A,B}}$ is directly
dependent on the security of $g^{ab}$, therefore for the sake of
simplicity we will consider $g^{ab}$ as the established key of the
exchange for the rest of this paper.

In the present article we will be concerned with the security of
this protocol. We will interpret security in a probabilistic manner
and will devise a statistical test that will ``assess'' the security
of the exchange in a given group.

In the cryptology literature there are two concepts of security --
the core security and the concept of semantic security which leads
to various security models. The semantic security and the related
concepts come under the name of ``provable security'' \cite[Section
2]{menezes-koblitz}. The core security of the Diffie-Hellman key
exchange protocol depends on the discrete logarithm problem, the
computational Diffie-Hellman problem and the decision Diffie-Hellman
problem. In this article we are concerned with the core security of
the exchange. We give a brief introduction to the discrete logarithm
problem and the computational Diffie-Hellman problem, for more on
these a reader can look at \cite[Section 5]{menezes-koblitz} or
\cite[Chapter 6]{stinson}.

\begin{assumption}[DL]
For a cyclic group $G$, generated by $g$, we are given $g$ and
$g^n$, $n\in\mathbb{N}$, the challenge is to compute $n$.
\end{assumption}

\begin{assumption}[CDH]
Given $g, g^a, g^b$ it is hard to compute $g^{ab}$.
\end{assumption}

Clearly, if these assumptions are not satisfied then $\mathcal C$,
an adversary\footnote{ There are various concepts of adversary in
cryptographic literature, the power and authority they have. In this
article we assume that our adversary is a passive eavesdropper.},
can gain access to the key $g^{ab}$. The relationship between these
two assumptions has been extensively studied. It is clear that the
CDH assumption will not be satisfied in a group where finding the
discrete logarithm solution is easy. In
\cite{Maurer+Wolf-RelaBetwBreaDiff:1999},
\cite{Boneh+Lipton-AlgoBlacFielthei:1996}, the authors show that in
several settings the validity of the CDH assumption and the hardness
of the discrete logarithm problem are in fact equivalent.

Unfortunately, the DL and the CDH assumptions are not enough to
ensure security of the Diffie-Hellman key exchange protocol. Even if
these assumptions are true, the eavesdropper $\mathcal C$ may still
be able to gain useful information about $g^{ab}$. For example, if
$\mathcal C$ can predict $90\%$ of the bits in $g^{ab}$ with high
probability then for all intents and purposes the key exchange
protocol is broken. Moreover, there exist protocols where the
knowledge of even one bit will break its security (some Casino
electronic games). With the current state of knowledge we cannot be
confident that assuming only CDH, a scenario like the one described
above does not exist (\cite{boneh98decision}).

\subsection{Our main contribution.}

Lemma \ref{lemma:DHIDDH} states that the security of the Diffie
Hellman exchange protocol is best studied from a statistical
perspective. We introduce a statistical treatment of this
particularly important problem in cryptography and it is our hope
that many more problems will be approached in a similar fashion.

We present novel methodologies to help asses the security of the
Diffie-Hellman key exchange protocol in a given group $(G,\cdot)$.
In Section \ref{sec:StatBackgr} we present the statistical criteria
we use as well as the relevance and connection with the security
assessment. Sections \ref{sec:testingDH-Indep} and
\ref{sec:testingDHI} present statistical tests to check the validity
of the statistical criteria presented in Section
\ref{sec:StatBackgr}. In particular, Subsection
\ref{subsec:permutTest} detail the use of the permutation testing
methodology to calculate concrete values for the probability of Type
I error of the tests. This section contains the important idea that
the method can be used to compare the security of the DH key
exchange protocol in two or more different groups. Furthermore, the
groups which we use to perform the comparison do not need to have
the same operational structure. Thus, in principle, it is possible
to compare the security of the exchange in finite groups generated
using elliptical curves versus the same order prime subgroups of
$\mathbb Z_n$, $n\in \mathbb N$. We do not pursue this direction in
the presented work.

Section \ref{sec:AppZp} applies the methodology we develop to some
examples where the security of the DH-exchange has been conjectured.
It is found that the results obtained strengthen the conjectured
hypotheses. Finally, in Section \ref{sec:conclusion} we present
general conclusions and directions of future research.

The treatment of the problem is based on the empirical distribution
of the key $g^{ab}$. We mention that a better approach from the
cryptographic perspective would be to look at the distribution of a
collection of bits in the binary expansion of $g^{ab}$. We believe
our methods could be extended and applied to this representation as
well.

\section{Statistical criteria to asses the security of the Diffie
 Hellman key exchange protocol.\label{sec:StatBackgr}}

In its most basic form described above the security of the
Diffie-Hellman key exchange protocol relies on an approximate
identification of the key $g^{ab}$ from the public information $g,\,
g^a,\, g^b$. In statistical terms there exist a clear concept that
answers the question of identification: statistical independence.
Therefore a \emph{sufficient} condition for the security of the DH
key exchange is:

\begin{assumption}[DH-Independence]
Given a cyclic group $G$ of order $N$, generated by $g$, let $a$ and
$b$ be chosen independently, uniformly at random from the set
$\{1,2,\dots,N\}$. Then the random variables $(g^a,g^b)$ and
$g^{ab}$ are independent.
\end{assumption}

For a given set $S$ we will use the notation $DU(S)$ to denote the
discrete uniform distribution on the elements of $S$. With this
notation $a$ and $b$ are independent random variables with the
$DU(\{1,2,\dots,N\})$ distribution.

Clearly this is a sufficient condition for the security of the
Diffie Hellman key exchange protocol. There is no information to be
gained about $g^{ab}$ from seeing $(g^a,g^b)$. Unfortunately, as
one's intuition may indicate, this assumption is rejected for any
finite group $G$ we have looked at. In the next section we construct
a statistical test for this assumption which will help introduce the
notations and the further testing procedures.

If the assumption presented above is not true, hope is not lost. The
DH-Independence assumption is a sufficient condition. In fact, in
the cryptographic literature this assumption is not even mentioned,
however a weaker necessary condition is presented:

\begin{assumption}[DDH]
Given $g, g^a, g^b$ and an element $z\in G$ it is hard to decide
whether or not $z=g^{ab}$.
\end{assumption}

In this form the DDH assumption constitutes a necessary condition
for the security of the Diffie-Hellman key exchange protocol.
Furthermore, \cite{joux:2003} construct groups based on elliptic
curves where the DDH assumption is not satisfied while the CDH and
the discrete logarithm problem are proven to be equivalent and hard.
This fact prompts the necessity to directly check the validity of
the DDH assumption for a given group.

The DDH assumption is assumed, either implicitly or explicitly in
many cryptographic systems and protocols. Applications include: the
many implementations of the DH key exchange itself (e.g.,
\cite{Diffie+OorschotETAL-Authauthexch:1992}), the El-Gamal
encryption scheme \cite{El-Gamal-Cryplogaoverfini:1984}, the
undeniable signatures algorithm
\cite{Chaum+Antwerpen-Undesign:1989}, Feldsman's verifiable secret
sharing protocol \cite{Feldman-PracScheNon-Inte:1987},
\cite{Pedersen-DistProvwithAppl:1991}, and many others; we point to
\cite{NaorDDHpseudorandom:1997} for a more detailed list.

Notice that the DDH assumption in the form presented above is a
little vague because of the use of the predicate, ``hard to
decide''. Surprisingly, attempts to make the DDH assumption explicit
were not made until late after its formulation in \cite{DH}. The
first ventures \cite[]{Boneh+Lipton-AlgoBlacFielthei:1996} use
standard cryptographic machinery (\cite{Yao-Theoappltrapfunc:1982,
Goldwasser+Micali-ProbEncr:1984}), to express the assumption in
terms of \emph{computational indistinguishability}. Put in this
traditional cryptographic form it was discovered quickly by
\cite{Stadler-PublVeriSecrShar:96} and independently
\cite{NaorDDHpseudorandom:1997} that if one assumes the existence of
a polynomial time probabilistic algorithm which distinguishes the
real key $g^{ab}$ from the other possible values even with a very
small probability\footnote{but not negligible. For the sake of
completeness we give here the whole definition. It is presented in
the footnote since it is not relevant to our approach at all.
Suppose that the group $G$ where the exchange takes place has order
$N$ and $n=\log_2 N$. It is said that a probabilistic algorithm
$\mathscr A$ decides on the right key with small (non-negligible)
probability if there exist a polynomial expression $p(\cdot)$ such
that for any $r\in G$: $$\left| \text{Prob}(\mathscr A \text{
outputs } g^{ab})- \text{Prob}(\mathscr A \text{ outputs } r)\right|
>\frac {1}{p(n)}.$$} (for all the possible inputs), then another
polynomial time algorithm can be constructed from the first which
will output $g^{ab}$ with a very large (almost one) probability. The
only requirement is that the size of the group is known, requirement
lessened by \cite{boneh98decision} which only requires finiteness of
the group.

All this evidence points toward a more specific definition based
entirely on the notion of \emph{statistical significance.} Indeed,
this fact materialized in a series of papers
\cite{canetti:1999,canetti:2000,FS,VNS}, which call this new form of
the assumption the Diffie Hellman Indistinguishability assumption
(DHI). We note that
\cite{Gennaro+KrawczykETAL-SecuHashDiffover:2004,joux:2003} use the
same form except it continues to call it DDH. We point the reader to
\cite{Hastad:1999} for a detailed discussion on the concept of
statistical significance versus computational significance; in the
context of pseudo-random number generation.

For our purposes of studying the security of the Diffie Hellman
exchange we will use the following assumption:
\begin{assumption}[DHI] Given $g, g^a, g^b$ the distribution of $g^{ab}$ is
indistinguishable from the Discrete Uniform distribution on the
elements of $G$ ($DU(G)$).
\end{assumption}
The notion of indistinguishability used here is the usual
statistical one. Two variables are indistinguishable if they have
essentially the same distribution, or put formally, $X_1$ and $X_2$
are indistinguishable if their distribution functions $F_i(x)=
P(X_i\leq x)$ with $i=1,2$ have the property:
$$F_1(x)=F_2(x), \text{ for all } x \in \mathbb R \setminus( A_1 \cup A_2), $$
where $A_1$, $A_2$ are the sets which contain the discontinuity
points of $F_1(\cdot)$, respectively $F_2(\cdot)$. Applied to our
specific case the distributions are discrete, therefore the
distribution functions $F_1$ and $F_2$ are just step functions with
jumps in a compact set in $\mathbb R$, thus using the right
continuity of the distribution functions, the usual definition
translates here in equality everywhere. We conclude that in our
context, indistinguishability means that the variables have the same
distribution.

This formulation is perfectly natural for a statistician who tries
to express the DDH formulation presented above. We note that our
version of the DHI assumption requires that the conditional
distribution $ (g^{ab} | g^a, g^b, g) $ is uniform while the
previous articles \cite{canetti:1999,canetti:2000,FS,VNS,
Gennaro+KrawczykETAL-SecuHashDiffover:2004,joux:2003} require that
the distribution of the triple $(g^a, g^b, g^{ab}|\,g)$ be Discrete
Uniform on the elements of $G\times G\times G=G^3$ ($DU(G^3)$).
Given an outcome $(x,y,z)$ we can write using the simple
multiplicative rule:
\begin{equation}
{{\mathbf P \left(g^a=x, g^b=y, g^{ab}=z|\, g \right)}}=\mathbf P
\left( g^{ab}=z |\, g^a=x, g^b=y, g \right) \mathbf P\left(g^a=x,
g^b=y |\, g \right)
\end{equation}
Under the original condition that $a$ and $b$ are
$DU(\{1,\ldots,N\})$ and using the fact that $g$ is a generator for
$G$ then the distribution of  $(g^a,g^b|\, g)$ is $DU(G^2)$, thus
the two formulations are perfectly equivalent.

It is known that in general \emph{statistical indistinguishability}
implies \emph{computational indistinguishability}, but the reverse
is not in general true, \cite[Section 3.2.2]{Goldreich-FounCryp:01}.
The following lemma states the same result in our specific case
using the assumptions presented in this section: DHI and DDH.

\begin{lemma}\label{lemma:DHIDDH}
In a group $G$ of order $N$, if the DHI assumption is true then the
DDH assumption is true as well.
\end{lemma}
\begin{proof}
Assume that DHI is true in $G$, then for given $g^a$, $g^b$, the
probability \\ $\mathbf{P}\left(g^{ab}=z|g^a,g^b\right)=1/N$ for any
$z\in G$. This is the hardest possible scenario in the DDH
assumption and hence we claim that DDH is satisfied.
\end{proof}
This lemma says that in any group $G$, DHI is a stronger\footnote{or
at least as strong} condition  than that of the DDH assumption. If
we look at the proof closely then we will see that the difference
between the DDH and the DHI consists in the fact that a measure of
hardness has been provided in the DDH assumption via the uniform
distribution.

\section{Testing for DH-Independence.\label{sec:testingDH-Indep}}

We give general definitions, then we go to our specific case.

Let $X$, $Y$, and $Z$ be three discrete random variables taking
values in the sets $\{x_1,x_2, \dots, x_n\}$,
$\{y_1,y_2,\dots,y_m\}$, $\{z_1,z_2,\dots,z_l\}$ respectively.
Denote with:
$$p(x_i,y_i,z_i) = \mathbf P \{X=x_i, Y=y_i, Z=z_i \},\quad \forall \, i,j,k$$
the joint probability function corresponding to $(X,Y,Z)$. With
usual notations we denote $p(y_j|x_i)$, $p(x_i,y_j|z_k)$, etc. the
conditional probability functions of $X|Y$, $(X,Y)|Z$, etc.
Furthermore, assume that for all $k\in\{1,2,\dots,l\}$ the marginal
distribution $p(z_k)=\mathbf P\{Z=z_k\}\neq 0$ to avoid
complications conditioning on a set of measure zero.

\begin{definition}[Entropy]
We define the {\it joint} and {\it conditional} {\bf measures of
uncertainty}.
\begin{eqnarray}
H\left( X,Y\right)& = &- \sum_{i=1}^{n}\sum_{j=1}^{m}\sum_{k=1}^{l}
p(x_i,y_j,z_k) \log p(x_i,y_j)  \label{eq:entropyxy} \\
H\left( X,Y|Z\right) & = & -
\sum_{i=1}^{n}\sum_{j=1}^{m}\sum_{k=1}^{l} p(x_i,y_j,z_k) \log
p(x_i,y_j|z_k)\label{eq:entropyxygivenz},
\end{eqnarray}
with the convention $0(-\infty)=0$.
\end{definition}

In the above definition we choose to work with the natural
logarithm, however any other basis will be equivalent for our
purpose due to the constant in the usual definition of the entropy
function (see \cite{S:1948}).

\begin{lemma}\label{lemma:entropy}
The following property holds for the above uncertainty measures:
$H(X,Y|Z)\leq H(X,Y)$ with equality if and only if $(X,Y)$ and $Z$
are independent.
\end{lemma}

The proof is an easy exercise in probability, the reader is directed
to \cite{S:1948} or \cite{R:1967} for more details.

Lemma \ref{lemma:entropy} gives a clear criterion for our first
test. More specifically: assume that the number of elements in $G$
is $N$, i.e. $|G|=N$. As an example $|\mathbb {Z}_p^\ast|=p-1$.

The plan is to apply the above lemma with $X=g^a$, $Y=g^b$ and
$Z=g^{ab}$. Since both participants in the Diffie-Hellman protocol
choose $a$ and $b$ at random and $g$ is the generator of $G$ we can
assume that $g^a$ and $g^b$ are independent and their distribution
is $DU(G)$. Thus, the distribution of $(g^a,g^b)$ is $DU(G\times
G)$. This in turn implies that $p(x_i,y_j)=1/N^2$ for all $i,j \in
\{1,\dots,N\}$, and thus the first entropy measure
(\ref{eq:entropyxy}) becomes:

\begin{equation}
H\left( g^a,g^b\right) = - \sum_{i,j,k} p(x_i,y_j,z_k) \log
\frac{1}{N^2} = 2 \log N
\end{equation}

At this point we can devise a test of the hypotheses:
\begin{equation}\label{test:independence}
\begin{cases}
H_0: & g^{ab} \text{ is independent of } (g^a,g^b)\\
H_a: & g^{ab} \text{ is NOT independent of } (g^a,g^b)
\end{cases}
\end{equation}
using Lemma \ref{lemma:entropy}. The test in
(\ref{test:independence}) is equivalent with:
\begin{equation}\label{test:entropy}
\begin{cases}
H_0: & H(g^a,g^b |\, g^{ab}) = 2\log N\\
H_a: & H(g^a,g^b |\, g^{ab}) < 2\log N
\end{cases}
\end{equation}

The question is: how do we proceed with this test? Since all the
distributions are finite, in theory at least, we could calculate
$p(x_i,y_j|z_k)$ for all the possible triples in $G\times G\times G
= G^3$. If we had these quantities it would be a simple matter to
calculate $H(g^a,g^b | g^{ab})$ according to
(\ref{eq:entropyxygivenz}).

Denote this value based on the whole set $G^3$ by $T_N$. The test
will then compare this value with $2\log N$. If equal then the
variables are independent and the DH-Independence assumption is
satisfied. If smaller then we could not prove independence of the
variables.

At this point let us make two important remarks.

\begin{remark}\label{remark:testisgeneral}
In the definition of the entropy functions (\ref{eq:entropyxy}) and
(\ref{eq:entropyxygivenz}) we did not use the structure of the group
$G$ in any way, only the relative frequency of the elements in the
group. This fact make the methods based on the entropy function well
suited for comparison between diverse groups. We will take advantage
of this feature later in this paper.
\end{remark}

\begin{remark}\label{remark:toomany}
In practice if we wish to calculate $T_N$ we have to calculate all
the possible values for $(g^a,g^b)$ and this will take longer than
an exhaustive search. Thus calculating $T_N$ is not practical,
instead we would have to estimate it. We will detail the estimation
in the next section.
\end{remark}

Alas, as we suspected from the beginning, implementing this first
test tells us that $(g^a,g^b)$ and $g^{ab}$ are not independent in
every group that we tried. For example in $\mathbb Z_p^\ast$ with
multiplication, calculating $T_N$ the entropy in
(\ref{eq:entropyxygivenz}) for $p\in \{1193,2131,11093\}$ will yield
values which are far apart from $2\log(p-1)$. In fact when looking
at the values obtained we see that they are close to $\log(p-1)$
thus the value of our first test increases with $p$. The closeness
of the test to $\log(p-1)$ is an interesting experimental fact. This
fact is investigated and explained by our second test presented in
the next section.

\section{Testing the DHI assumption \label{sec:testingDHI}}

If the DH-Independence assumption is satisfied in a given group $G$,
then we could stop and decide that we found a perfect group for the
Diffie Hellman key exchange. However, the experimental procedures
and our intuition point out that the DH-Independence assumption is
never satisfied in any \emph{finite} group $G$. The next task is to
obtain a statistical testing procedure to verify the validity of the
DHI assumption in a given group $G$. The idea is to use the entropy
function \eqref{eq:entropyxygivenz} in the sense of Kullback-Leibler
divergence \cite{kullback:1951} as a measure of departure from the
entropy calculated under the hypothesis of Uniform distribution.
Specifically, using earlier notation, we wish to construct a
statistical test that will check the validity of the following
hypotheses:

\begin{equation}\label{test:uniformity}
\begin{cases}
H_0: & \text{The distribution of } g^{ab}|\, (g^a,g^b) \text{ is } DU(G)\\
H_a: & \text{The distribution of } g^{ab}|\, (g^a,g^b) \text{ is NOT
}  DU(G)
\end{cases}
\end{equation}

Let us denote the elements of $G$ as $\{g_1,g_2,\dots,g_N\}$.
Suppose we can look at all the possible triples $(g^a,g^b,g^{ab})$
when $a,b\in \{1,2,\dots,N\}$ take all the possible values. Clearly,
there are $N^2$ such possible triples and assuming that $a$ and $b$
are chosen at random, each such triple will have probability
$1/N^2$. The last element in the triple $g^{ab}$ will get mapped
into $N$ possible values (the elements of $G$). Thus, some values in
$G$ will be repeated. For an element $g_k\in G$ denote $m_k$ the
number of times $g_k$ appears in the place of $g^{ab}$ among all the
$N^2$ triples. We have then $\sum_{k} m_k=N^2$. For any pair
$(g^a,g^b)$ that corresponds to $g^{ab}=g_k$ we can then calculate
the conditional probability as:
$$ p(g^a=g_i,g^b=g_j |\,
g^{ab}=g_k)=\frac{1}{m_k}\mathbf{1}_{A}(g_i,g_j,g_k),$$ where $A$ is
the set of all possible $N^2$ tuples $\left(g^a,g^b,g^{ab}\right)$,
and we have used the notation $\mathbf 1_A(x)$ to denote the
indicator function of the set $A \subset \Omega$, i.e., $\mathbf 1_A
: \Omega \rightarrow \{0,1\}$ is given by:
$$\mathbf 1_A(x)=\begin{cases}1 & \quad \text{ if } x\in A \\
0 & \quad \text{ if } x\not\in A \end{cases}$$ We can continue:
\begin{align}
  H\left( g^a , g^b| g^{ab}\right) & = - \sum_{i=1}^{N}\sum_{j=1}^{N}\sum_{k=1}^{N}
p(g_i,g_j,g_k) \log p(g_i,g_j|g_k)\nonumber\\
& = - \sum_{k=1}^{N}\sum_{i,j=1}^{N}\frac{1}{N^2}\log
\frac{1}{m_k}\,\mathbf 1_{A}(g_i,g_j,g_k)\nonumber\\
& = -\sum_{k=1}^{N}
\frac{m_k}{N^2}\log\frac{1}{m_k}.\label{eq:entropyours}
\end{align}

Under the null hypothesis $H_0$, the distribution of $(g^{ab}|\,
g^a, g^b)$ is uniform, therefore we should have the $m_k$
multiplicities equal. This automatically implies that $m_k=N$ for
all $k$'s and then the entropy function in \eqref{eq:entropyours}
is:
\[ H\left( g^a , g^b| g^{ab}\right)= -\sum_{k=1}^{N}
\frac{1}{N}\log\frac{1}{N} = \log N\]

The testing statistics is:
\begin{equation}\label{eq:DHItest}
T_N = H\left( g^a , g^b|\, g^{ab}\right)-\log N = \sum_{k=1}^{N}
\frac{m_k}{N^2}\log m_k -\log N.
\end{equation}
This test is based on the whole set of values in $G^2$. Accordingly,
if the value of the test equals zero then the null hypothesis $H_0$
is true, any other value of the test will support the alternative
hypothesis. We summarize this result in the following:

\begin{lemma}[Testing Procedure]
With the previous notations if $T_N = 0$ then the DHI assumption is
satisfied in a given group $G$. \end{lemma}

Both remarks \ref{remark:testisgeneral} and \ref{remark:toomany}
certainly apply for this testing procedure as well. In particular,
remark \ref{remark:toomany} means that we have to find procedures to
estimate $T_N$ instead of calculating it. This will introduce
distributions and we detail the approach next.

\subsection{The permutation test approach.\label{subsec:permutTest}}
Assume that we can obtain a sample of $n$ pairs $\{(a_i,b_i)\}_{i
\in \{1,2,\dots,n\} }$ from $\{1,2,\dots,N\}\times \{1,2,\dots,N\}$.
For each pair in the sample we can calculate the triple $(g^{a_i},
g^{b_i}, g^{a_i b_i})$. Let $A_n$ be the set of all the triplets in
the sample.

Using \eqref{eq:entropyours} we can calculate an estimate of
$H\left( g^a , g^b| g^{ab}\right)$ using:
\begin{align}\label{eq:probestimates}
  \hat p_n(g_i, g_j, g_k) & =\frac{k_{ijk}}{n} \mathbf 1_{A_n}(g_i, g_j, g_k) ,\\
 \hat p_n(g_i, g_j|\, g_k) & = \frac{k_{ijk}}{m_k} \mathbf 1_{A_n}(g_i, g_j, g_k),\nonumber
\end{align}
where once again $m_k$ denotes the multiplicity of $g_k$, but in the
given sample of $n$ observations. We took into account the
possibilities of obtaining repeated observations in the sample by
multiplying with the factor $k_{ijk}$; which represents the number
of times we see the same observation $(g_i, g_j, g_k)$ in our
sample.

The test statistic is:
\begin{equation}\label{eq:TestStatOurs}
T_n = - \sum_{i=1}^{n}\sum_{j=1}^{n}\sum_{k=1}^{n} \hat
p_n(g_i,g_j,g_k) \log \hat p_n(g_i,g_j|g_k) -\log n
\end{equation}

All that is left, is to investigate the distribution of $T_n$ under
the null hypothesis $H_0$. Under the null hypothesis the $m_k$'s are
the multiplicities of $g_k$'s in a sample of size $n$ drawn from the
set $$\{g_1,\dots,g_1,g_2,\dots,g_2,\dots,g_N,\dots,g_N\}$$ where
each element in the group $G$ are repeated $N$ times.

Let us denote $M_1,M_2,\dots,M_N$ the multiplicities of the elements
$\{g_1,g_2,\dots,g_N\}$ in a sample of size $n$. It is not hard to
show that the joint probability distribution of $(M_1,\dots,M_N)$ is
the so called \emph{multivariate hypergeometric distribution}:
\[\mathbf P \left( M_1=m_1,\dots,M_N=m_N\right) =
 \frac{{N \choose {m_1}}{N \choose {m_2}}\dots{N \choose {m_N}}}
 {{{N^2} \choose n}}\]

The test statistic under $H_0$ is:

\begin{equation}\label{eq:TestStatunderH0}
\widehat{T}_n= \sum_{k=1}^{N} \frac{M_k}{n}\log{M_k} -\log n.
\end{equation}

If we would be able to calculate the distribution of $T_n$ knowing
that $(M_1,M_2,\dots,M_N)$ are multivariate hypergeometric then we
would be in position to reach the conclusion of the test of
uniformity \eqref{test:uniformity} by calculating the p-value of the
test statistic \eqref{eq:TestStatOurs} using this distribution.

Finding the distribution of the test statistic under $H_0$
\eqref{eq:TestStatunderH0} is however not an easy task. This is the
reason we propose the use of permutation testing for which knowledge
of this distribution is not necessary.

The permutation testing procedure generates samples
$(M_1,M_2,\dots,M_N)$ from the Multivariate hypergeometric
distribution. For each sample, it calculates the corresponding value
of the test statistic under the null hypothesis as in
\eqref{eq:TestStatunderH0}. These values are obtained from the
assuminption that $H_0$ is true; this allow us to calculate the
empirical distribution of our sample statistic $T_n$ under the null
hypothesis. The $p$-value of our test is given by the proportion of
values as extreme or more than the one calculated in
\eqref{eq:TestStatOurs} using the group $G$.

A small $p$-value is an evidence against the null hypothesis in
\eqref{test:uniformity}, that the sample comes from a uniform
distribution. We summarize the procedure bellow:

\newcounter{ocount}
{\bf Testing procedure to determine validity of DHI for a group $G$}
\begin{list}{(\roman{ocount})}{\usecounter{ocount}}
\item We take a sample of size $n$ and we calculate the
test statistic as in \eqref{eq:TestStatOurs}.
\item We generate many test statistic values under the hypothesis
$H_0$ is true using \eqref{eq:TestStatunderH0}, then construct their
empirical distribution.
\item We calculate the $p$-value of the test as the proportion of
values in the empirical distribution found in (ii) lower than the
test value found using $G$ in (i).
\item If the $p$-value is small we reject the DHI assumption. If the
$p$-value is big we did not find evidence that the DHI is not
satisfied in the given group $G$.
\end{list}

\subsection{How to compare two or more groups?}

We will note at this point that the absolute value of the test
$|T_N|$ and its estimate $|T_n|$ represent a measure of departure
from the Discrete Uniform distribution. The bigger the estimate the
further is the distance from the uniform distribution and the weaker
is the validity of the DHI assumption. Remark
\ref{remark:testisgeneral} also tells us that the nature of the
group operation is irrelevant for the testing procedure. Therefore,
we can use the test as a tool to compare the strength of the
Diffie-Hellman key exchange protocol in two or more groups. To be
able to do so we need the order of the groups compared to be similar
and, more importantly, the sample size on the basis of which we
calculate the permutation test to be the same. We take advantage of
the ability to compare different groups in the next section.

\section{Testing the DHI assumption in ${\mathbb Z_p^\ast}$}\label{sec:AppZp}

We are going to check the efficiency of the testing procedure for
the most useful finite groups, those included in ${\mathbb
Z_p^\ast}$ with the multiplicative operation. We present the
following examples as a way for checking the validity of the testing
procedure.

\begin{example}[A group where the DDH assumption does not
hold.]\label{example:DHnotsecure} Consider $G={\mathbb Z_p^*}$ with
$p$ prime. It is known that computing Legendre symbol in this group
gives a distinguisher against DDH
(\cite{Gennaro+KrawczykETAL-SecuHashDiffover:2004}).
\end{example}

\begin{example}[A group where the DDH assumption is conjectured to
hold]\label{example:DHsecure} We currently do not know any DDH
distinguisher for a prime order subgroup of ${\mathbb Z_p^*}$.
Therefore, given $p$ and $q$ prime with $g$ divisor of $p-1$ it is
conjectured that in a subgroup of order $q$ of ${\mathbb Z_p^*}$ the
DDH assumption holds.
\end{example}

We start with a given group $G$ and using the test presented in the
previous section we will test for the validity of the DHI assumption
in that group $G$. This should provide a strong indication towards
the security of the Diffie-Hellman key exchange protocol in that
group.

\subsection{The rate of convergence of the testing procedure}

\begin{figure}[h]
\centering
\includegraphics[width=11cm]{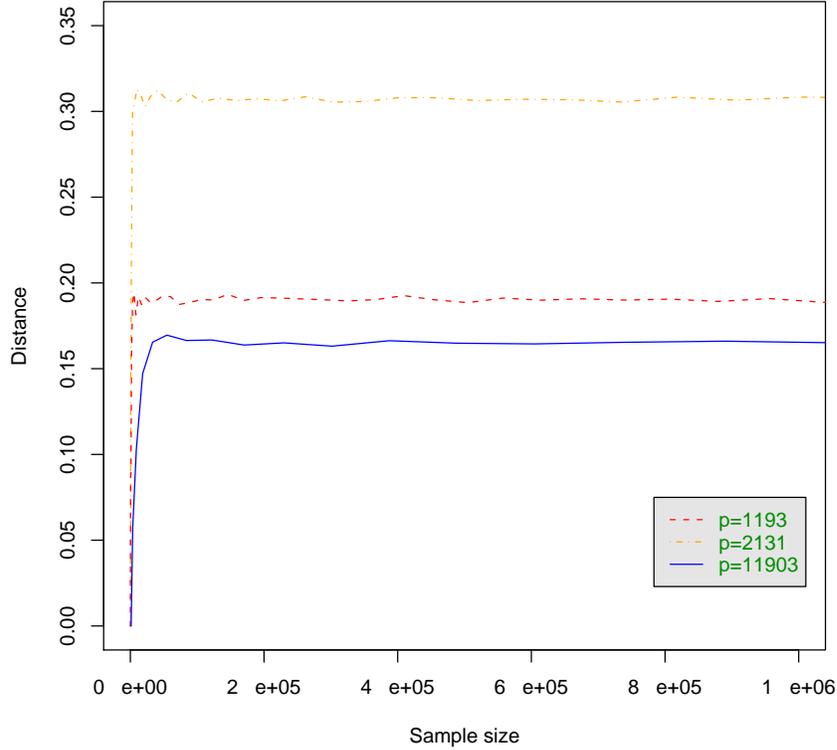}
 \caption{Comparison of the test values for different sample sizes
  and ${\mathbb Z_p^*}$'s}\label{fig:distances}
\end{figure}

The first thing we investigate is the rate of convergence for our
test. To do this we need to calculate the true value of $T_N$ and
thus we have to look at small groups.

For space consideration we are presenting only results obtained for
$p=1193$ in Table \ref{table:results1193} in the Appendix. The
sample sizes are presented in the first column of the table and the
corresponding sample entropy value $T_n$ in column two. Column 3
presents the proportion of values lower than $T_n$ -- an entry equal
to $1$ corresponds to a $p$-value 0 of the test. The fourth value in
the table represents the distance from $T_n$ to the center of the
distribution of entropy values calculated under $H_0$. Finally, the
last value represents the ratio of the distance in column four, to
the distance from the sample entropy $T_n$ to the furthest away
point in the distribution. It is an indication on how many standard
deviations away $T_n$ is from the distribution.

There are two remarkable features of these values -- one, we see
that the test rejects the null hypothesis that the distribution of
$(g^{ab}| \, g^a,g^b)$ is uniform on the elements of $G$; the other
remarkable feature is that we determine this fact based on a sample
of $354$ values or about one third of the value of $N=1192$. In the
second place if we wish to determine the actual entropy distance
from the two distributions -- a feature that will be useful when
comparing two or more groups; we can see that starting with a sample
size of $n=3304$ (or about $3$ times $N$) we start to obtain
accurate results.

To illustrate better the rate of convergence for some other groups
we plot in Figure \ref{fig:distances} the evolution of the test
values with the size of the sample. This figure suggest that to get
a good estimate for $T_N$ the sample size will depend on the size of
the group, for example we need a larger sample size for ${\mathbb
Z_{11903}^\ast}$ than we need for ${\mathbb Z_{1193}^\ast}$. In
addition, the same figure points out another interesting fact.

Following example \ref{example:DHnotsecure} we know that ${\mathbb
Z_p^\ast}$ is not secure. It is also conjectured that some groups
are more secure than others. Looking at the problem from that
perspective, for which groups are more easily broken using the
Legendre symbol, it is also assumed that by increasing the size of
the group one can make the group more secure.

We can see from the figure that the second assertion is not true.
Just increasing the size of the group does not make it more secure.
Remembering that a smaller relative distance corresponds to
closeness to the Discrete uniform distribution on the elements of
$G$, we see from the Figure \ref{fig:distances} that while ${\mathbb
Z_{11903}^*}$ the largest group is the most secure of the three, the
situation between the other two groups is not what we would have
expected looking at the size of the group alone. Even though
${\mathbb Z_{2131}^*}$ is the larger group (almost twice the size),
it is also less secure from the DHI assumption perspective than
${\mathbb Z_{1193}^*}$. This indicate that the choice of the group
$G$ rather than the size of it is essential for the security of the
Diffie-Hellman key exchange protocol.

\subsection{Comparison of the DHI assumption across groups.}

Next we wished to give an indication of groups that are more secure
than others. It is known that considering only the Legendre symbol
criterion the safest groups among ${\mathbb Z_p^*}$ are the ones
obtained when $p$ is a safe prime i.e., of the form $p=2q+1$ where
$q$ is another prime \cite{menezes:1996}.

\begin{figure}[h]
\centering
\includegraphics[width=11cm]{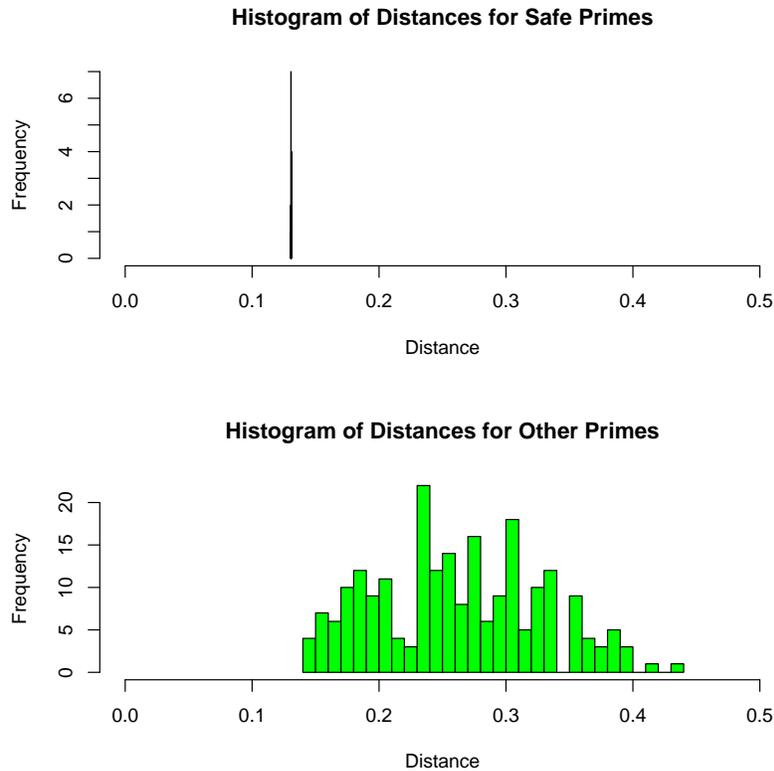}
\caption{Histogram of all the test values for ${\mathbb Z_p^*}$ with
$2000<p<4000$. Values closer to zero represent safer groups for DH
exchange.\label{fig:histSafePrimes.p2000}}
\end{figure}

\begin{figure}[h]
\centering
\includegraphics[width=11cm]{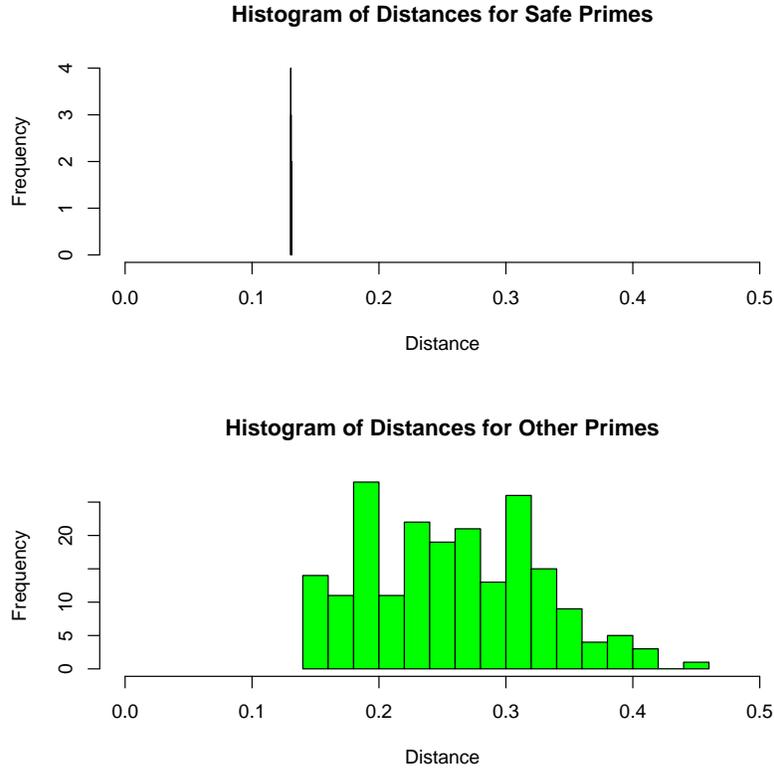}
\caption{Histograms of test values obtained for ${\mathbb Z_p^*}$
with $9000<p<11000$. Values closer to zero represent safer groups
for DH exchange. \label{fig:histSafePrimes.p9000}}
\end{figure}

\begin{figure}[h]
{\centering
\includegraphics[width=11cm]{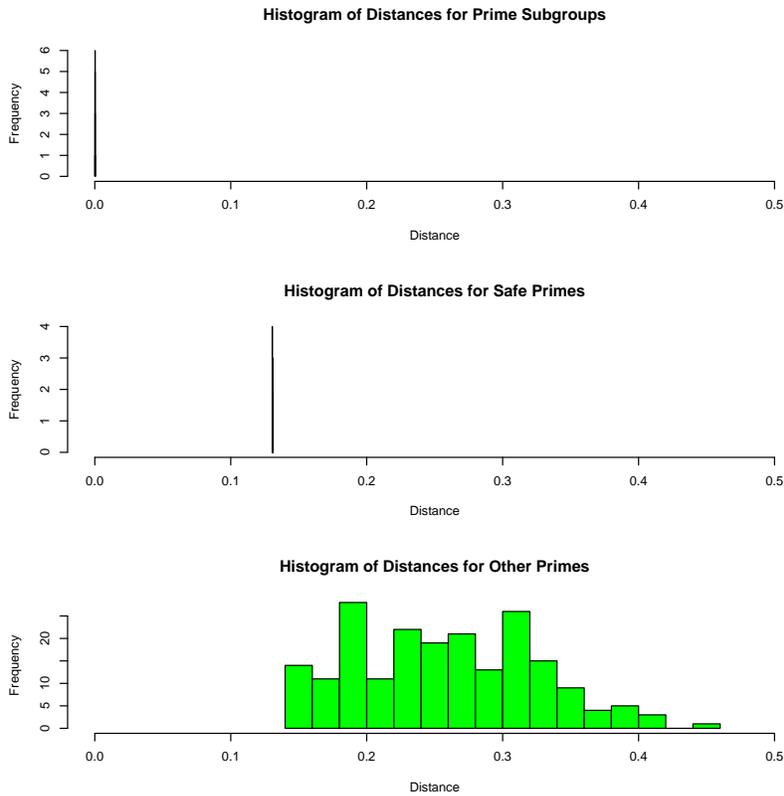}
\caption{Comparing values of the test for different type of groups
when $9000<p<11000$. On top, we plot values for prime subgroups of
${\mathbb Z_p^*}$ when $p$ is a safe prime. Middle, we plot values
for safe prime ${\mathbb Z_p^*}$'s. On bottom, values for all the
other groups ${\mathbb Z_p^*}$ in the range given. Values closer to
zero represent better groups for DH exchange.
\label{fig:histSubgroupsSafe} }}
\end{figure}

\begin{figure}[h]
\centering
\includegraphics[width=11cm]{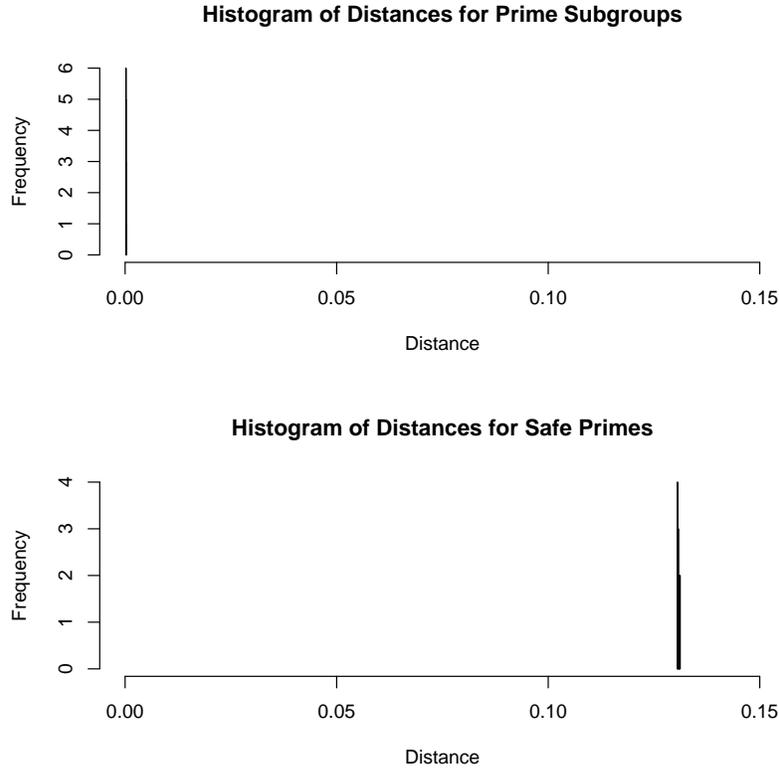}
\caption{A more detailed comparison of the previous image (Fig.
\ref{fig:histSubgroupsSafe}). We are comparing the prime subgroups
with the corresponding safe groups. Values closer to zero represent
safer groups for DH exchange. \label{fig:histSafevsSubgroups}}
\end{figure}

\begin{figure}[h]
\centering
\includegraphics[width=11cm]{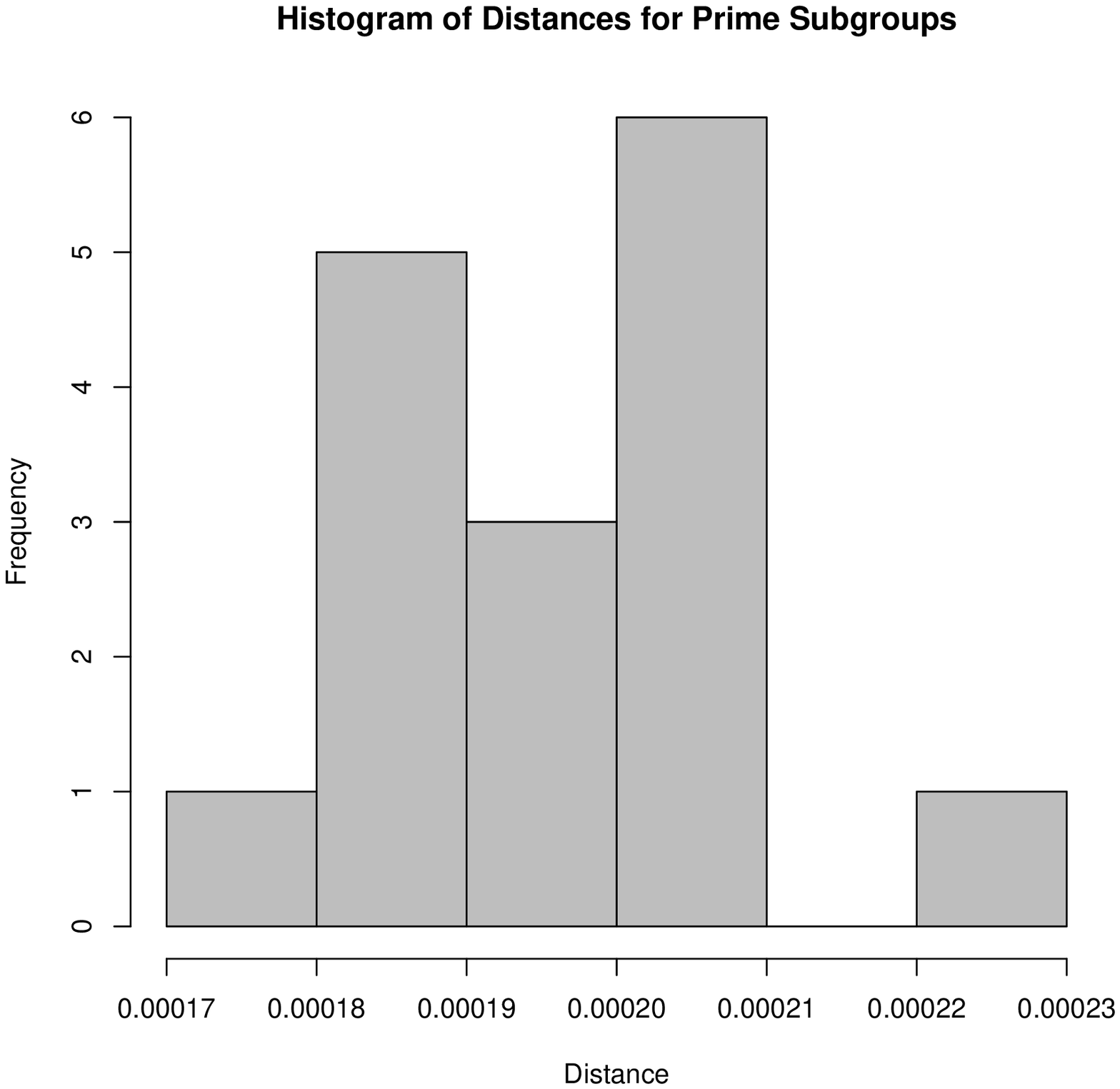}
\caption{A blowup of the histogram of the values for the prime
subgroups in the safe primes. Note the values are close to zero but
not equal to zero. \label{fig:histSubgroupSafePrimes}}
\end{figure}

We wished to test this theory for a large set of ${\mathbb Z_p^*}$
groups with varying $p$'s. We looked at all primes between 2000 and
4000, and again for primes between 9000 and 11000. The reason for
the two separate segments of primes is that we expect some sort of
consistency between them. We show the distribution of the test
values for these groups separated into safe and not safe primes in
Figures \ref{fig:histSafePrimes.p2000} and
\ref{fig:histSafePrimes.p9000}.

First, we notice that the behavior of primes in the range 2000 to
4000 is very similar with the primes for the higher range 9000 to
11000. Second, in both ranges we see the same conclusion applies,
the safe prime groups are more secure than any other groups.
However, the test estimate obtained for each of the safe prime
groups is significantly different from zero therefore there is no
safe group in the ranges given for which the DHI assumption is
verified. This seem to confirm the assertion in the Example
\ref{example:DHnotsecure}.

Next, we will look to Example \ref{example:DHsecure}. We will use
our test for the prime subgroups of each of the safe primes in the
range 9000 to 11000. More specifically, we look at each ${\mathbb
Z_p^*}$ with $p$ a safe prime, and we construct the prime subgroup
of order $q$ in each such group. Then we test the DHI assumption in
each subgroup thus constructed. The values obtained for the
distances are plotted in the upper histogram of Figure
\ref{fig:histSubgroupsSafe}. We mention that the behavior of the
test values for primes between 2000 and 4000 was very similar, for
space consideration we omit the corresponding plot. All the values
are obtained using the same sample size $n=8\times 10^6$. The reason
for this particular value is that while the groups themselves are in
the range $9000$ to $11000$, the subgroups are of order $4500$ to
$5500$.

It is remarkable to see that these subgroups are clearly safer for
the DH exchange than any other groups plotted in the picture. The
results seem to confirm the conjecture in the Example
\ref{example:DHsecure}. However, the actual test of uniformity was
rejected, but we needed a very large sample size almost equal to the
maximum value $N^2$.

For a better comparison we plotted in Figure
\ref{fig:histSafevsSubgroups} on page
\pageref{fig:histSafevsSubgroups} only the histogram of the values
obtained for the prime subgroups of the $Z_p^\ast$ with $p$ a safe
prime (top) and the histogram of the values obtained for the
$Z_p^\ast$ groups, $p$ a safe prime between $9000$ and $11000$
(bottom).

It is remarkable the closeness of these values to each other
considering that the order of the group varies between $9000$ and
$11000$ a $20\%$ variation in size. This is an encouraging fact,
which suggests that for even larger $p$'s we will see the same sort
of consistency in the values. This will imply that groups with the
same operational structure will have similar behavior from the point
of view of the Diffie-Hellman security. However, there is a
variation in the values as illustrated in the Figure
\ref{fig:histSubgroupsSafe} on page \pageref{fig:histSubgroupsSafe}
which represent the histogram of the values obtained for the prime
subgroup of $Z_p^\ast$ groups, with $p$ a safe prime varying between
$9000$ and $11000$.

\section{Conclusion and future work.\label{sec:conclusion}}

In this article we present a novel statistical testing procedure to
help assess the security of the Diffie-Hellman key exchange
protocol. The methods presented are quite general and to our
knowledge represent the first systematic pure applied statistical
approach to a cryptographic problem. The article is intended to open
a way for methods coming from statistical world to the cryptographic
domain. We do not claim to solve the security of the Diffie-Hellman
exchange protocol. What we have presented are primarily sufficient
conditions for the security. We also presented a way to compare the
strength of these conditions in different groups. In Section
\ref{sec:AppZp} we show that among the groups we looked at, only the
prime subgroups of a large group are close to fulfilling the
conditions considered.

An obvious lack in our results is a statistical analysis for very
large primes. Typically the groups used in cryptography are of the
order at least $2^{1024}$. The use of our testing procedure,
ad-literam as presented in section \ref{sec:testingDHI} prevents us
from such an analysis, however currently we are investigating
directions of circumventing the permutation testing approach. One
direction is to approximate the distribution of the test in
\eqref{eq:TestStatunderH0} with a multinomial distribution, then use
a multivariate normal distribution for a second approximation. This
should allow us to calculate the p-value of the test directly
without the need of the permutation testing. Another direction is to
put together outcomes into coarser groups and look at the
distribution of these groups of outcomes. This idea is similar in
result with the approach of \cite{canetti:1999, banks:2006}, and
should allow us to speed up the procedure in order to apply it to
much larger groups. It will also allow us to look at the
distribution of the binary representation of prime subgroups of a
large group, and extend the methodology to finite groups defined
using elliptical curves.

\section{Appendix}
We present the actual values obtained in $\mathbb Z_p$ when $p=1193$
in Table \ref{table:results1193}.

\begin{center}
\newcommand\T{\rule{0pt}{2.6ex}}
\newcommand\B{\rule[-1.2ex]{0pt}{0pt}}
{\small
\begin{longtable}{c|cccc}
\caption[Results $\mathbb Z_{1193}^\ast$]{Results for $\mathbb
Z_{1193}^\ast$\label{table:results1193}}\\
  \hline
Sample size $n$ \T \B & Sample entropy value& $p$-value& Distance to center
& Relative distance\\
  \hline
\endfirsthead
\hline \multicolumn{5}{|l|}%
{{\bfseries \tablename\ \thetable{} -- Continued from previous
page}\T\B} \\ \hline Sample size $n$ \T \B & Sample entropy value&
$p$-value& Distance to center
& Relative distance \\
  \hline
\endhead

\hline \multicolumn{5}{|r|}{{Continued on next page}} \\ \hline
\endfoot

\hline \hline
\endlastfoot

59 \T &0.046993&0.556&0&0\\118&0.105734&0.904&0&0\\
354&0.280115&1&0.0867205869841293&0.602176619941004\\
885&0.532425&1&0.088729342259792&0.599662439145006\\
1829&0.96382&1&0.158395336513686&0.758210206038124\\
3304&1.40654&1&0.187918729961140&0.890429900342397\\
5428&1.82531&1&0.194266177572582&0.935988768705612\\
8319&2.19741&1&0.181355529391107&0.952411456549936\\
12095&2.55884&1&0.192004761885750&0.966980890255525\\
16874&2.87286&1&0.187337211259202&0.979981576295216\\
22774&3.1674&1&0.191522958630831&0.981517705830948\\
29913&3.43077&1&0.188465586935031&0.98875796589123\\
38409&3.67754&1&0.189706561218385&0.989433516631953\\
48380&3.90781&1&0.192416008075197&0.99165124060656\\
59944&4.11938&1&0.19204137533093&0.99478104337756\\
73219&4.31302&1&0.187468331653386&0.994755817256502\\
88323&4.50025&1&0.188526799093478&0.996560990251708\\
105374&4.67745&1&0.19031690687346&0.996854516374416\\
124490&4.84304&1&0.190069869416784&0.997275184401866\\
145789&5.00357&1&0.193382446593161&0.997349334367475\\
169389&5.14947&1&0.189808592541566&0.997931831476642\\
195408&5.29352&1&0.191440055573543&0.998590870156603\\
223964&5.42925&1&0.191148605525655&0.998411134116565\\
255175&5.55893&1&0.190698438462096&0.998845341861062\\
289159&5.68315&1&0.190158376706143&0.998956167303273\\
326034&5.80232&1&0.189542921352024&0.99921427003345\\
365918&5.91821&1&0.190224685898690&0.999099568272852\\
408929&6.03153&1&0.192585342387550&0.99931590521473\\
455185&6.13611&1&0.190147697214377&0.999413667702933\\
504804&6.2378&1&0.188502041891033&0.999534588253575\\
557904&6.34038&1&0.191174900210519&0.999622308526143\\
614603&6.43583&1&0.189935678716927&0.999625957882701\\
675019&6.53032&1&0.190747522920258&0.99976657291958\\
739270&6.62041&1&0.189993001424682&0.99970573975723\\
807474&6.70913&1&0.190533231900746&0.999813379231167\\
879749&6.79349&1&0.189228458019329&0.999814900993474\\
956213&6.87841&1&0.190858660137478&0.99984125150255\\
1036984&6.95729&1&0.188695781636576&0.999914186731075\\
1122180&7.03801&1&0.190502664400493&0.999926167521206\\
1211919&7.11461&1&0.190209994699298&0.99994757464427\\
1306319 \B &7.18871&1&0.189337283972479&0.999979948549426\\

\end{longtable}}
\end{center}


\bibliographystyle{chicago}
\bibliography{ddh}

\end{document}